\documentclass[english, a4paper, 11pt]{article}
\pdfoutput=1
\usepackage[T1]{fontenc}
\usepackage[utf8]{inputenc}
\usepackage[margin=3cm]{geometry}
\usepackage[
    unicode=true,
    pdfstartview=FitV,
    colorlinks=true,
    citecolor=Turquoise,
    linkcolor=Turquoise,
    urlcolor=MFCB,
    linktoc=page,
    hyperindex=true,
    pdfcreator={},
]{hyperref}
\usepackage[english]{babel}
\usepackage{amsmath,amsfonts,amssymb}
\usepackage{amsthm}
\usepackage[english]{cleveref} 
\crefname{figure}{Figure}{Figures}
\Crefname{figure}{Figure}{Figures}
\crefname{section}{Section}{Sections}
\crefname{Rq}{Remark}{Remarks}
\crefname{Def}{Definition}{Definitions}
\crefname{Ex}{Example}{Examples}
\crefname{CEx}{Counter-example}{Counter-example}
\theoremstyle{definition}
\newtheorem{Def}{Definition}
\numberwithin{Def}{section}
\newtheorem{Ex}[Def]{Example}

\newtheorem{Rq}[Def]{Remark}
\numberwithin{equation}{section}
\usepackage{lmodern}
\usepackage{bm}
\usepackage{eulervm}
\usepackage[x11names]{xcolor} 
\usepackage[pdftex]{graphicx}
\usepackage{subcaption} 
\usepackage{tikz} 
\definecolor{MCB}{cmyk}{0,0.03,0.08,0.26} 
\definecolor{MFCB}{cmyk}{0,0.06,0.20,0.6} 
\colorlet{Turquoise}{DeepSkyBlue4}
\colorlet{TurquoiseClair}{DeepSkyBlue4!20}
\colorlet{Orange}{DarkOrange3!85}
\usepackage[frozencache=true,cachedir=minted-cache]{minted} 
\usepackage{setspace} 
\usepackage{enumitem}
\setlist{nosep}

\title{Drawing Diestel-Leader graphs in 3D}
\author{Amandine Escalier}
\date{\today}
\hypersetup{ 
  pdfauthor={Amandine Escalier},
  pdftitle={Drawing Diestel-Leader graphs in 3D},
  pdfsubject={TikZ code to draw Diestel-Leader graphs in 3D},
  pdfkeywords={Diestel-Leader, TikZ, graphs, 3D}
}
\newcommand{\deffont}[1]{\textbf{\textcolor{DeepSkyBlue4!75!black}{#1}}} 
\newcommand{\frakh}{\mathfrak{h}}
\newcommand{\bfT}{\mathbf{T}}
\newcommand{\bN}{\mathbb{N}}
\newcommand{\bZ}{\mathbb{Z}}

\begin{document}

\maketitle
\begin{abstract}
  In this short note we give some code to represent Diestel-Leader graphs in 3D.
  The code is written in TikZ. 
\end{abstract}

\begin{figure}[htbp]
  \centering
  \begin{subfigure}[b]{0.47\textwidth}
    \centering
    \includegraphics[width=\textwidth]{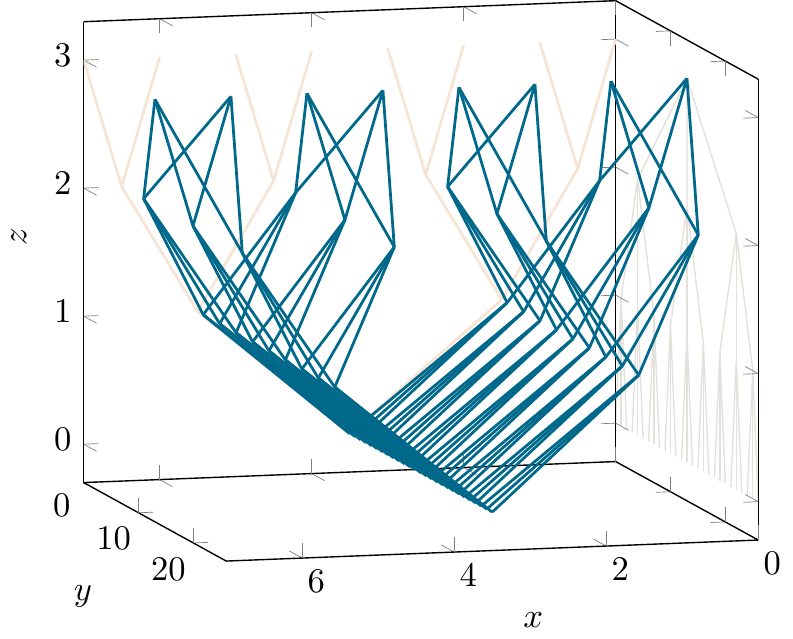}
    \label{subfig:DL32scale31}
  \end{subfigure}\hfill
  \begin{subfigure}[b]{0.47\textwidth}
    \centering
    \includegraphics[width=\textwidth]{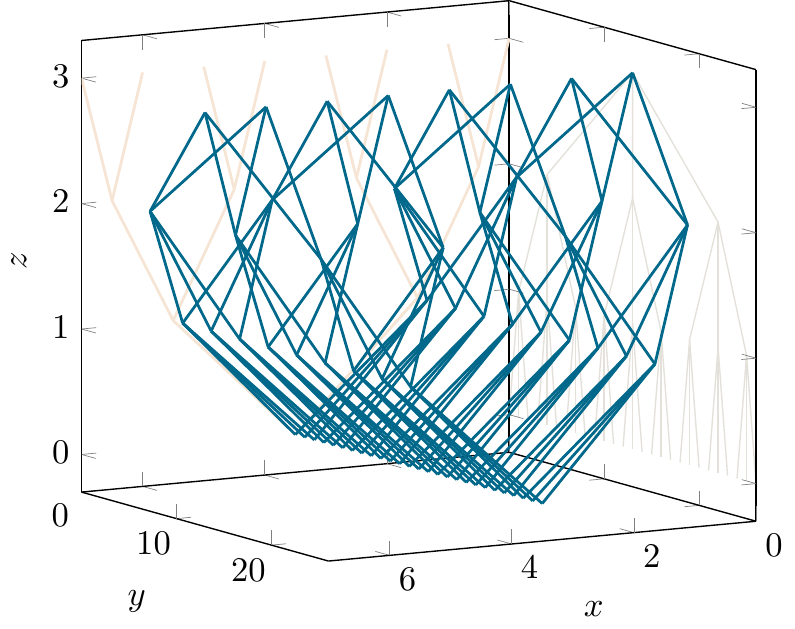}
    \label{subfig:DL32scale32}
  \end{subfigure}
  \caption[DL(3,2)]{Two views of the Diestel-Leader graph $DL(3,2)$}
  \label{fig:DL32}
\end{figure}

The history of Diestel-Leader graphs takes its root in the following question asked
by Woess \cite{Woess}: is every connected locally ﬁnite vertex-transitive graph
quasi-isometric to some Cayley graph? In the hope of answering no to this
question, Diestel and Leader \cite{DL} defined 
what we now call \emph{Diestel-Leader graphs}. However it was only later, in the
famous papers of Eskin, Fisher an Whyte \cite{EFW0,EFW1} that it was showed that
some of the aforementioned graphs are not quasi-isometric to any Cayley graph.

In this note we give a code to draw these graphs in 3D, using TikZ. Readers only
interested by producing an illustration of some $DL(p,q)$ can jump to the last
pages of this article (or the end of the .tex file) and copy-paste the
code given in \cref{Code}, write the wanted values of $p$ and $q$ (line 29) and
then compile. Readers wishing to change the code can rely on the description
made in \cref{subsec:Comments}. We start this note with a short reminder of the
definition of Diestel-Leader graphs.

\newpage
\section{Diestel-Leader graphs}
We recall here the definition of Diestel-Leader graphs.
We refer to \cite[Section 2]{Woess}
for more details.

\subsection{Tree and horocycles}\label{subsec:DefHorocycles}
Let $q\geq 2$ and denote by $\bfT  =T_q$ the homogenous tree of degree $q+1$. Denote by
$d$ the usual graph distance on $\bfT$ fixing to $1$ the length of an edge.

A \deffont{geodesic ray} is an infinite sequence $(v_n)_{n\in \bN}$ of vertices
of $\bfT$ such that $d(v_i,v_j)=|i-j|$, for all $i,j\in \bN$. We say that two
rays are equivalent if ther symmetric difference\footnote{Recall that the
  symmetric difference of two sets $A$ and $B$ is defined by $|A\triangle
  B|=(A\backslash B)\cup (B\backslash A)$} is finite.
We call \deffont{end} of $T_q$ an equivalence class of rays in $\bfT$ and denote
by $\partial \bfT$ the space of ends of $\bfT$. 

Let $\hat{\bfT}:=\partial \bfT \cup \bfT$. For any elements $x,y\in
\hat{\bfT}$ there is a unique geodesic in $\hat{T}$, denoted by \deffont{$\overline{xy}$},
that connects $x$ and $y$.

\begin{figure}[htbp]
  \centering
    \includegraphics[width=0.85\textwidth]{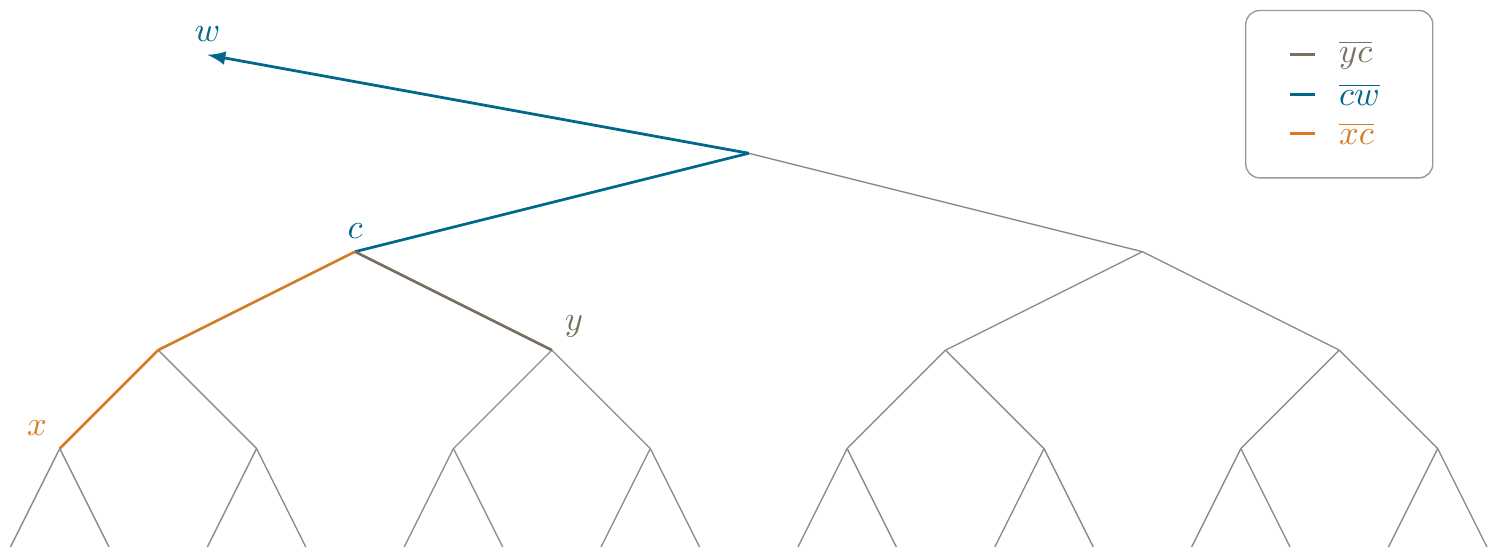}
  \caption[Example of confluent]{Example of the confluent of two points in $T_2$}
  \label{fig:Confluent}
\end{figure}

Now fix an end $w\in \partial \bfT$, the \deffont{confluent} of two elements
$x,y\in \hat{\bfT}\backslash {w}$ with respect to $w$, denoted by $x\curlywedge y$, is
defined as the element $c=x\curlywedge y$ such that $\overline{xw}\cap
\overline{yw}=\overline{cw}$, that is to say: the confluent is 
the point where the two geodesics $\overline{xw}$ and $\overline{yw}$ towards
$w$ meet (see \cref{fig:Confluent}).

Now, fix a root vertex $o\in \bfT$. We define below a \emph{Busemann
  function}, which will allow us to endow our tree with some height notion.

\begin{Def}\label{Def:Busemann}
  Let $w\in \partial T_q$. The \deffont{Busemann function} with respect to $w$
  is the map $\frakh:\bfT \rightarrow \bZ$ defined by
  \begin{equation*}
    \frakh(x)=d(x,x\curlywedge o)-d(o,x\curlywedge o).
  \end{equation*}
\end{Def}

\begin{Ex}
 Let’s turn back to \cref{fig:Confluent} and let $o=y$ be the root. Then for $x$
 represented in the figure, we have $x\curlywedge o=c$ and thus
 $d(x,x\curlywedge o)=2$ and $d(o,x\curlywedge o)=1$. Therefore $\frakh(x)=1$.
\end{Ex}

\begin{figure}[htbp]
  \centering
  \includegraphics[width=0.95\textwidth]{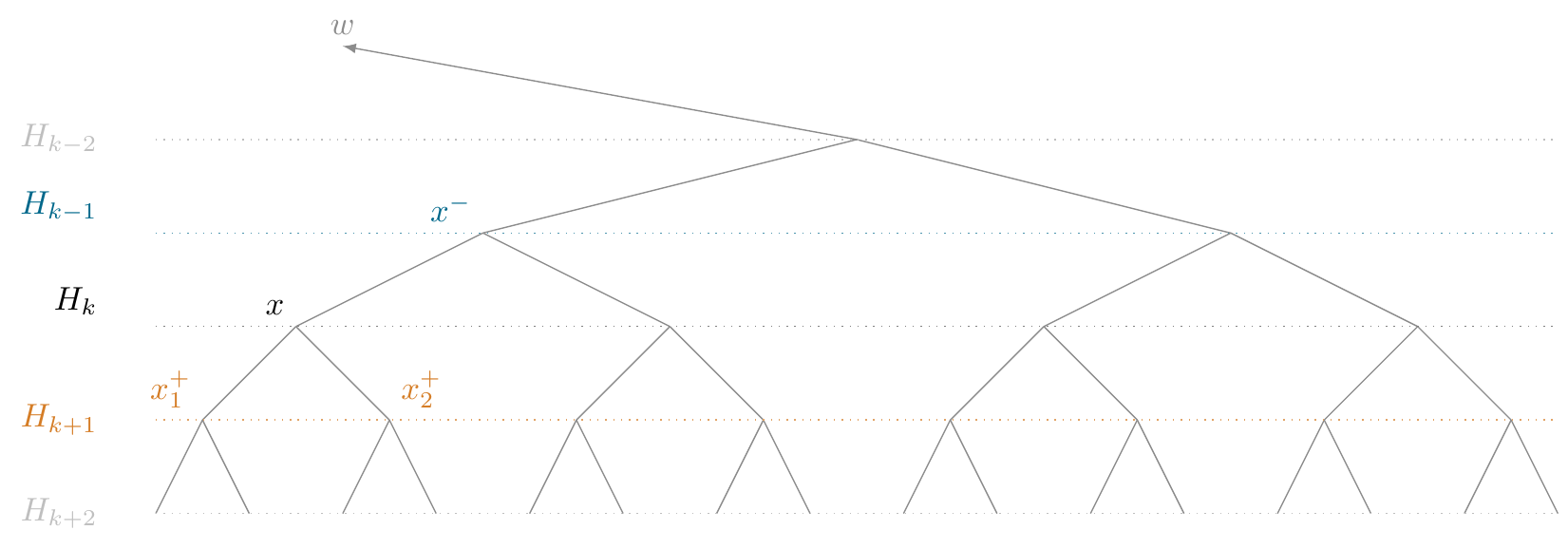}
  \caption{Horocycles, predecessor, successors}
  \label{fig:Horocycles}
\end{figure}

\begin{Def}
  Let $w\in \partial T_q$ and $k\in \bN$. The \deffont{horocycle} with respect
  to $w$, denoted by $H_k$, is the set $H_k=\left\{ x\in \bfT:\frakh(x)=k
  \right\}$.
\end{Def}
We refer to \cref{fig:Horocycles} for an illustration.
Note that every horocycle in $\hat{T}$ is infinite. Every vertex $x$ in a
horocycle $H_k$ has exactly one neighbour $x^-$ (called the
\deffont{predecessor}\label{Def:predecessor}) in $H_{k-1}$ and exactly $q$ neighbours
(called \deffont{successors}) in $H_{k+1}$ (see 
\cref{fig:Horocycles} for an illustration).

\subsection{Diestel-Leader graphs}\label{subsec:DefDLpq}
Now fix $p,q \geq 2$ and consider $T_p$ and $T_q$ with respective roots $o_p$
and $o_q$, and respective reference ends $w_p$ and $w_q$.
\begin{Def}
  The Diestel-Leader graph $DL(p,q)$ is the graph whith set of vertices
  \begin{equation*}
    V\left(DL(p,q)\right)=\big\{ (x,y)\in T_p\times T_q : \frakh(x)=-\frakh(y)\big\}.
  \end{equation*}
  and where there is an edge between two elements $(x_1,y_1)$ and $(x_2,y_2)$,
  if and only if $(x_1,x_2)$ is an edge in $T_p$ and $(y_1,y_2)$ is an edge in
  $T_q$. 
\end{Def}

\begin{figure}[htbp]
  \centering
  \begin{subfigure}[b]{0.48\textwidth}
    \centering
    \includegraphics[width=\textwidth]{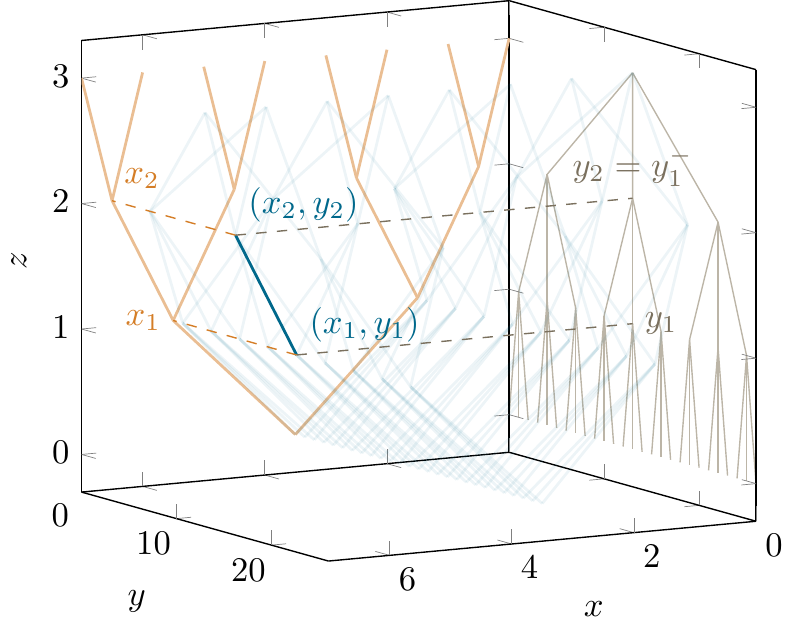}
    \caption{When $x_2$ is a child of $x_1$ and $y_2=y^{-}_1$}
    \label{subfig:Edge1}
  \end{subfigure}\hfill
  \begin{subfigure}[b]{0.48\textwidth}
    \centering
    \includegraphics[width=\textwidth]{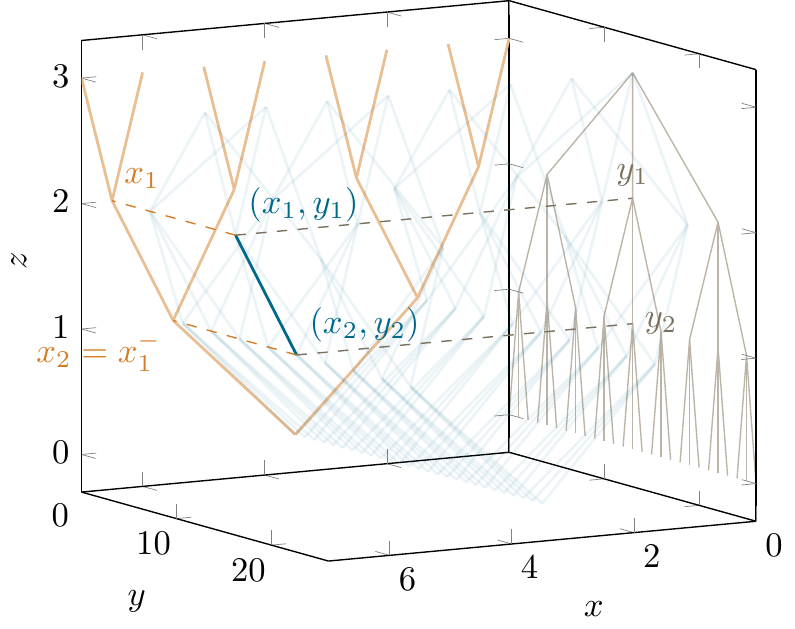}
    \caption{When $y_2$ is a child of $y_1$ and $x_2=x^{-}_1$}
    \label{subfig:Edge2}
  \end{subfigure}
  \caption[Two examples of edges]{Drawing an edge in $DL(p,q)$: two cases}
  \label{fig:Edges}
\end{figure}

If there is an edge between two vertices $(x_1,y_1)$ and
$(x_2,y_2)$ in $DL(p,q)$, then remark that either
\begin{itemize}
\item $x_2$ is one of the $p$ childs of $x_1$ in $\bfT_p$ and in this case
  $y_2$ is the only predecessor of $y_1$ in $\bfT_q$, namely $y_2=y^{-}_1$ (see \cref{subfig:Edge1});
\item or $x_2=x^{-}_1$ is the unique predecessor of $x_1$ in $\bfT_p$
  and in this case $y_2$ is one of the $q$ childs of $y_1$ in $\bfT_q$ (see
  \cref{subfig:Edge2}).
\end{itemize}

\begin{figure}[htbp]
  \centering
  \begin{subfigure}[b]{0.48\textwidth}
    \centering
    \includegraphics[width=\textwidth]{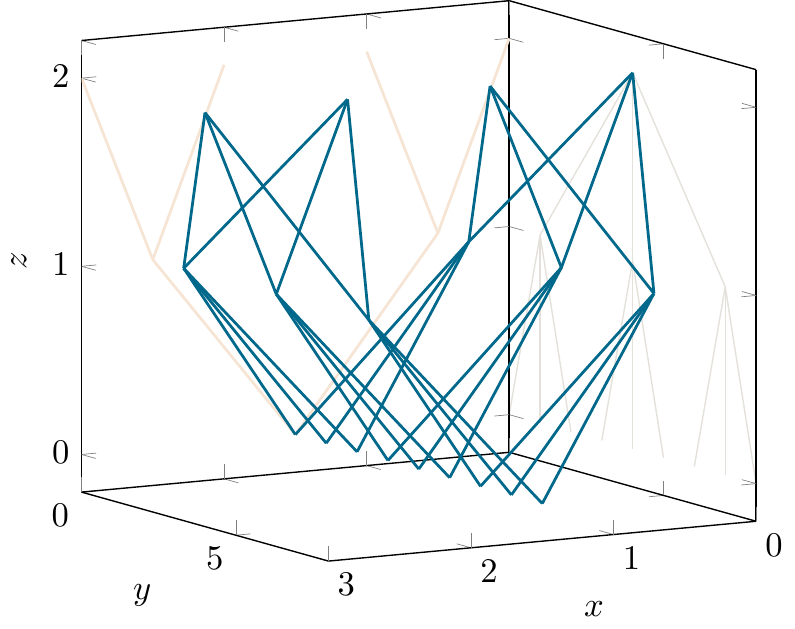}
    \caption{$DL(3,2)$}
    \label{subfig:DL32Scale2}
  \end{subfigure}\hfill
  \begin{subfigure}[b]{0.48\textwidth}
    \centering
    \includegraphics[width=\textwidth]{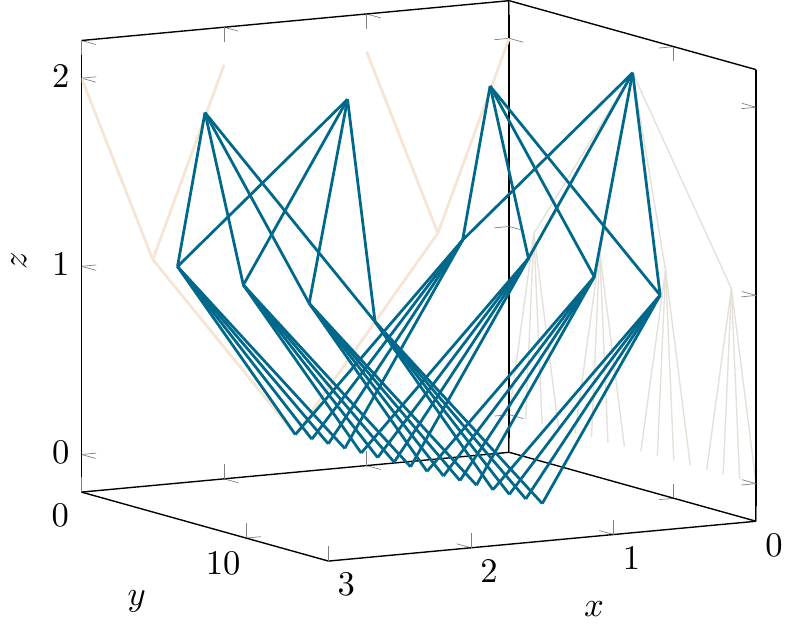}
    \caption{$DL(4,2)$}
    \label{subfig:DL42scale2}
  \end{subfigure}
  \caption[Two examples]{Two different Diestel-Leader graphs, represented at the same
  scale}
  \label{fig:DL32andDL42}
\end{figure}

We represent these two cases in \cref{fig:Edges}. The edge is drawn in
\textcolor{DeepSkyBlue4}{blue}, the $p$-regular tree in
\textcolor{Orange}{orange} and the $q$-regular tree in
\textcolor{MFCB}{brown}. In light blue is represented the Diestel-Leader
graph, we refer to \cref{fig:DL32andDL42} for a
drawing of $DL(p,q)$ itself. Note that in \cref{fig:DL32andDL42}, the
corresponding $T_p$ and $T_q$ are drawn on the planes $x=0$ and $y=0$
respectively.

\begin{Rq}
  When $p=q$ the graph $DL(p,q)$ is a Cayley graph of the Lamplighter group
  $\bZ/p\bZ \wr \bZ$.
\end{Rq}
\section{The code}
The complete code is given in \cpageref{Code}, and written in TikZ. We start by
some comments on how the coordinates were computed and how the loops work.

\subsection{Comments on the code}\label{subsec:Comments}

\paragraph{Variables} 
The main variables are 
\mintinline{latex}{\p, \q} (line 29) and \mintinline{latex}{\layers} (line
35). The first two variables correspond to the number of childs in the first and
the second tree respectively, that is to say the $p$ and $q$ in $DL(p,q)$.
Finally, the height of the graph is parametrized by \mintinline{latex}{\layers}:
there are \mintinline{latex}{\layers}$+1$ different heights of vertices
appearing on the graph.

With the values given in the example in \cref{Code}, we produce the left image in
\cref{fig:DL32}.

Finally changing the values in \mintinline{latex}{view={165}{10}} (line $17$)
will change the point of view (see for example \cref{fig:DL32}). For more
details on how to parametrize the view, see the \emph{pgfplots} manual
\cite[Section 4.11.1, page 311]{PGFPlots}.

\paragraph{Structure of the code}
The code is composed of one main loop cut in two parts: the first one that draws
the tree on the plane $y=0$, and the second one that draws simultaneously the
tree located on the plane $x=0$ and the Diestel-Leader graph.

In this description “height” will stand for the $z$-coordinates in the picture. 

\begin{description}
\item[Lines 58 to 72] We start by drawing the regular tree of degree $p$
represented on the plane $y=0$, in orange in the pictures. This tree grows from
the bottom to the top.

At height \textcolor{DeepSkyBlue4!60}{$h$} there are
$p^{\textcolor{DeepSkyBlue4!60}{h}}$ vertices to draw and the
space between two consecutive nodes at this height is 
$\textcolor{black!50}{\mathrm{pspace}}=p^{\text{layers-\textcolor{DeepSkyBlue4!60}{h}}}$
(in particular at the top, when $\textcolor{DeepSkyBlue4!60}{h}=layers$, the
space between two nodes is equal to $1$).

Moreover the first horizontal node at height $\textcolor{DeepSkyBlue4!60}{h}$ has for
horizontal coordinates the middle of
$[0,p^{\mathrm{layers}-\textcolor{DeepSkyBlue4!60}{h}}-1]$ that is to say
\begin{equation*}
  p^{\mathrm{layers}-\textcolor{DeepSkyBlue4!60}{h}}/2-0.5=\textcolor{black!50}{\mathrm{pspace}}/2-0.5,
\end{equation*}
and therefore, if $\textcolor{DeepSkyBlue4!60}{h}\neq layers$, the first of its
child (at height $\textcolor{DeepSkyBlue4}{n}=\textcolor{DeepSkyBlue4!60}{h}+1$) 
has coordinates the middle of
$[0,p^{\mathrm{layers}-\textcolor{DeepSkyBlue4}{n}}-1]$, namely
$p^{\mathrm{layers}-\textcolor{DeepSkyBlue4!60}{h}}/(2p)-0.5%
=\textcolor{black!50}{\mathrm{pspace}}/(2p)-0.5$.

\begin{figure}[htbp]
  \centering
  \includegraphics[width=\textwidth]{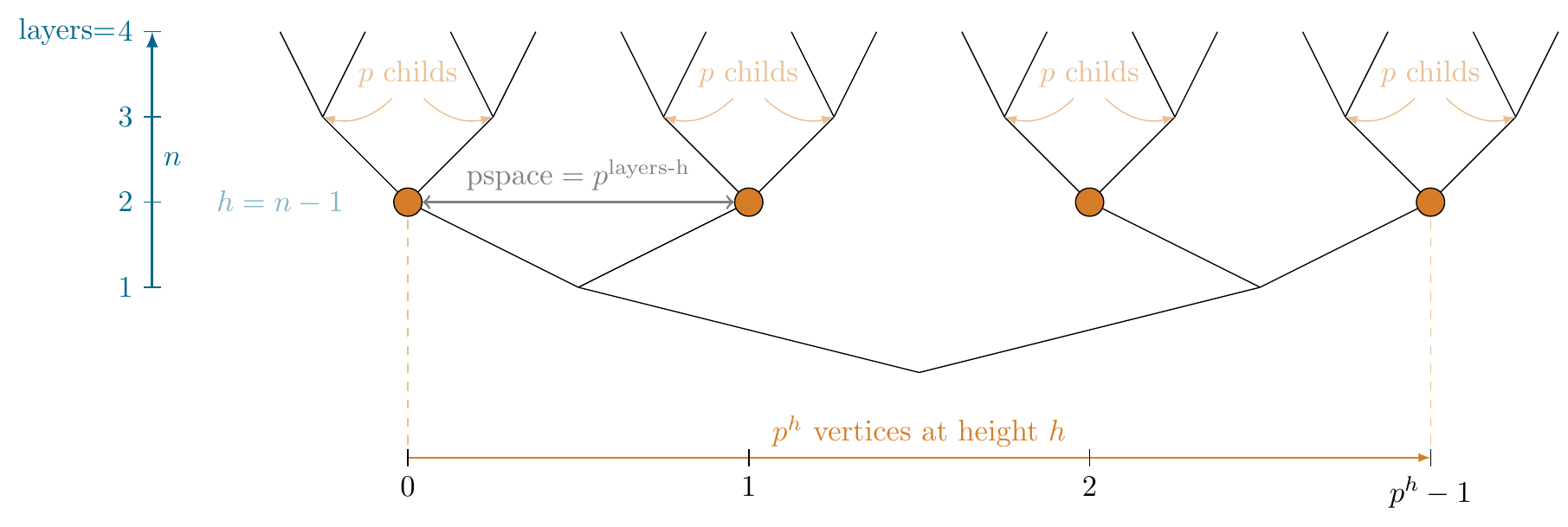}
  \caption{Ranges of the different loop variables for the first tree}
  \label{fig:Expli1}
\end{figure}

Now fix \textcolor{DeepSkyBlue4}{$n \in \{1,\ldots,\text{layers}\}$}, the loop
starting at line $58$ goes over all the nodes at height
$\textcolor{DeepSkyBlue4!60}{h}=\textcolor{DeepSkyBlue4}{n}-1$ and draws the
edges between a node at height $n-1$ and its $q$ childs at height $n$. To do so, it
fixes the horizontal position \textcolor{Orange}{$k\in \{1,\ldots,p^{n-1}-1\}$}
of the chosen node at height $n-1$. This node thus has to be shifted to the
right by $\textcolor{Orange!50!black}{pshift}
=\textcolor{Orange}{k}*\textcolor{black!50}{\mathrm{pspace}}$ 
and hence drawn at the coordinates

\begin{equation*}
  (\textcolor{black!50}{\mathrm{pspace}}-0.5 
  +\textcolor{Orange}{k}*\textcolor{black!50}{\mathrm{pspace}},0,n-1)
  = (\textcolor{black!50}{\mathrm{pspace}}-0.5
  +\textcolor{Orange!50!black}{pshift},0,n-1).
\end{equation*}
The childs are then located at height $\textcolor{DeepSkyBlue4}{n}$. Similarly
as above, the first child to be drawn has for horizontal coordinates the middle of
$[0,q^{\mathrm{layers}-\textcolor{DeepSkyBlue4}{n}-1}-1]$ and the
“$\backslash\mathrm{child}$-th” is shifted
to right by $\mathrm{child}*\textcolor{black!50}{\mathrm{pspace}}$. Hence the
coordinates of the “$\backslash\mathrm{child}$-th” child is given by 
\begin{equation*}
  (\textcolor{black!50}{\mathrm{pspace}}/q-0.5
  +\textcolor{Orange!50!black}{pshift}
  +\mathrm{child}*\textcolor{black!50}{\mathrm{pspace}},0,n-1).
\end{equation*}

\item[Lines 79 to 116] This loop draws the brown tree on the plane $x=0$ and
the Diestel-Leader graph \emph{at the same time}.
\begin{description}
\item[Lines 79 to 95] The variable $k$ at line $79$ goes through the
  $\textcolor{MFCB}{qrev}=q^{\mathrm{layers}-\textcolor{DeepSkyBlue4}{n}}$
  vertices at height $n$. Denote by $v$ the $k$-th vertex at this height. The
  loop \mintinline{latex}{for \child in {0,...,\q-1}} goes through the $q$ childs
  of $v$ and draw the edge between $v$ and the $\backslash \mathrm{child}$-th
  child at height $\textcolor{DeepSkyBlue4}{n}-1$. The coordinates of the nodes
  are computed in the same 
  way as in the previous loop. (Remark that the heights in the current loop are
  swapped: the parent is at height $n$ and the childs are at height
  $\textcolor{DeepSkyBlue4}{n}-1$ whereas in the previous one, the childs were
  at height $\textcolor{DeepSkyBlue4}{n}$ and the parent at height
  $\textcolor{DeepSkyBlue4}{n}-1$.)
\item[Lines 99 to 114] \emph{Inside} the loop
  \mintinline{latex}{for \child in {0,...,\q-1}}
  (starting line 81) lies also the part that draws the Diestel-Leader. Remember
  that $v$ is the current vertex at height $\textcolor{DeepSkyBlue4}{n}$ in the
  $q$-regular tree $\bfT_q$.
  
  The loop \mintinline{latex}{for \kk in {0,...,\pnm-1}} goes through all the
  vertices in the orange tree, drawn at height $n$. It will correspond to all
  the couple $(uv)$ such that $u\in \bfT_p$ is drawn at the same height
  (the $z$-coordinates) than $v$. To a given child $c$ of $v$ in $\bfT_p$,
  corresponds then the edge in the Diestel-Leader linking $(u,v)$ to
  $(u^{-},c)$, where we recall that $u^{-}$ is the unique predecesor of $u$ in
  $\bfT_p$ (see \cpageref{Def:predecessor}). Hence the point $(u,v)$ will be
  drawn at height $n$ and will have the same $x$-coordinates than $u$ and same
  $y$-coordinates than $v$, namely
  
  \begin{equation*}
    \begin{cases}
      x&=\textcolor{black!50}{\mathrm{pspace}}/(2p)-0.5
      +\textcolor{Orange!50!black}{\mathrm{pshiftprime}}
      +childprime*\textcolor{black!50}{\mathrm{pspace}}/p,\\
      y&=q^n/2-0.5+qshift,\\
      z&=n
    \end{cases}
  \end{equation*}
  Similarly, the node $(u^{-},c)$ will be drawn at height $n-1$ and will have the
  same $x$-coordinates than $u^{-}$ and same $y$-coordinates than $c$, namely
  \begin{equation*}
    \begin{cases}
      x&=\textcolor{black!50}{\mathrm{pspace}}/2-0.5 
      +\textcolor{Orange!50!black}{\mathrm{pshiftprime}},\\
      y&=q^n/2-0.5+\mathrm{qshift}+\mathrm{child}*q^n/q,\\
      z&=n-1
    \end{cases}
  \end{equation*}
\end{description}
\end{description}

\subsection{The code}\label{Code}
\begin{minted}[linenos]{latex}
  \documentclass[border=0mm]{standalone}
  \usepackage[x11names]{xcolor}
  \usepackage{tikz}
  \usetikzlibrary{calc,math,backgrounds}
  \usepackage{ifthen}
  \usepackage{pgfplots}
  \pgfplotsset{compat=1.18}
  % -----------------------------------------------------------------
  % ---------------------------- Colors (Colorblind friendly)
  % -----------------------------------------------------------------
  \definecolor{MFCB}{cmyk}{0,0.06,0.20,0.6} 
  \colorlet{Orange}{DarkOrange3!85}
  % + DeepSkyBlue4
  % -----------------------------------------------------------------
  \begin{document}
  \begin{tikzpicture}%[scale=0.5] %Uncomment to change the scale
    \begin{axis}[
      view={165}{10}, % Change the point of view.
      % view={150}{10}, % Change the point of view.
      xlabel=$x$,
      zlabel=$z$,
      ylabel=$y$,
      ]
      \tikzmath{
        % ------------------------------------
        % ------------------------------------
        % ---- PARAMETERS p AND q OF DL(p,q)
        % ------------------------------------
        % ------------------------------------
        \p=2; \q=3;
        % ------------------------------------
        % ------------------------------------
        % ----- HEIGHT PARAMETER
        % ------------------------------------
        % ------------------------------------
        \layers=3;% If layer=3 then there are
        %4 heights appearing, numbered 0 1 2 3
        % ------------------------------------
        % ------------------------------------
        for \n in {1,...,\layers}{% Vertical, n stands for the height
          % -----------------------------------------
          % Stored variables
          % -----------------------------------------
          % % For the q regular tree
          \qrev=pow(\q,\layers-\n);% Stores the value 2^(L-n)
          \qn=pow(\q,\n);%Stores the value q^n
          % For the p regular tree drawn from bottom to top
          \pspace=pow(\p,\layers+1-\n); % (space between two nodes
          % at height layers-(n-1))
          \pnm=pow(\p,\n-1);% Stores p^(layers-nprime)=p^(n-1)
          % --------------------------------------------------------
          % Regular tree of degree p drawn on the plane y=0
          % Drawn starting from the bottom to the top
          % --------------------------------------------------------
          for \k in {0,...,\pnm-1}{%Horizontal
            \pshift=\k*\pspace;% Horizontal shift: k*p^n
            for \child in {0,...,\p-1}{% Child of the considered node
              {
                % draw an edge
                \addplot3[Orange!20,thick] coordinates%
                {%The vertex at the top (ie. the child) at height n
                  (\pspace/(2*\p)-0.5+\pshift+\child*\pspace/\p,0,\n)
                  % The vertex below (ie. the parent) at height n-1
                  (\pspace/2-0.5+\pshift,0,\n-1)};
              }; % End of edges drawing
            }; % End of the loop “for \child in”
          };% End of the loop “for \k in ”
          % =========================================================
          % Regular tree of degree q drawn on the plane x=0
          % ---------------- AND -----------------------------------
          % The Diestel-Leader graph
          % =========================================================
          for \k in {0,...,\qrev-1}{%Horizontal
            \qshift=\k*\qn;% Horizontal shift: k*q^n
            for \child in {0,...,\q-1}{%Goes through the q childs
              %in the q-regular tree
              %
              % The Tree drawn on the plane x=0 (the light brown one)
              % 
              { %Drawing the edge
                \begin{scope}[on background layer]
                  \addplot3[MFCB!20,thick] coordinates%
                  % vertex at height n
                  {(0,\qn/2-0.5+\qshift,\n)
                    % child at height n-1
                    (0,\qn/(2*\q)-0.5+\qshift+\child*\qn/\q,\n-1)};
                \end{scope}
              };%End of drawing for the edge of the tree
              % 
              % The Diestel-Leader
              % 
              for \kk in {0,...,\pnm-1}{
                \pshiftprime=\kk*\pspace;% horizontal shift
                for \childprime in {0,...,\p-1}{%Goes through the p
                  % childs in the p-regular tree
                  {%
                    % draw a blue edge of the Diest-Leader 
                    \addplot3[DeepSkyBlue4,thick] coordinates%
                    {%The vertex at the top
                      (\pspace/(2*\p)-0.5+\pshiftprime%
                      +\childprime*\pspace/\p,
                      \qn/2-0.5+\qshift,\n)
                      % The vertex at height n-1
                      (\pspace/2-0.5+\pshiftprime,\qn/(2*\q)%
                      -0.5+\qshift+\child*\qn/\q,\n-1)};
                  };%End drawing
                };% End if the loop “for childprime in”
              };% End if the loop “for kk in”
            }; % End of the loop “for \child in”
          };% End of the loop “for k in”
        };% End of the loop “for n in”
      }%End Tikzmath
    \end{axis}
  \end{tikzpicture}
\end{document}
\end{minted}

\newpage
\bibliographystyle{alpha}
\bibliography{DLpq.bib}
\vfill
\textbf{Acknowledgments} 
We thank Tom Ferragut, Giles Gardam and Théo Laurent for their useful remarks
and comments.
\bigskip

\noindent \textbf{Fundings} The athor is funded by the Deutsche
Forschungsgemeinschaft (DFG, German Research Foundation) – Project-ID
427320536 – SFB 1442, as well as under Germany’s 
Excellence Strategy EXC 2044 –390685587, Mathematics Münster:
Dynamics–Geometry– Structure.
\vfill
\noindent\textbf{\textcolor{DeepSkyBlue4}{Amandine Escalier}}\\
Mathematisches Institut,\\
Fachbereich Mathematik und Informatik der Universität Münster,\\
Orléans-Ring 12,\\
48149 Münster,\\
Germany
\end{document}